\documentclass[11pt]{amsart}
\usepackage[T1]{fontenc}
\usepackage[utf8]{inputenc}
\usepackage[francais,english]{babel}
\usepackage{frbib}
\usepackage{amsmath} 
\usepackage{graphicx} 
\usepackage{geometry}           
\usepackage{subfig} 
\usepackage{amssymb}
\usepackage{amsthm}
\usepackage{epstopdf}
\usepackage{enumerate}
\usepackage{hyperref}
\usepackage{pinlabel}

\hypersetup{
backref=true, 
pagebackref=true,
hyperindex=true, 
breaklinks=true, 
urlcolor= blue, 
linkcolor= blue, 
bookmarks=true, 
bookmarksopen=true, 
pdftitle={Propriétés  combinatoires du  bord d'un groupe hyperbolique - Séminaire TSG Grenoble}, 
pdfauthor={Antoine Clais}, 
pdfsubject={} 
}

\newtheorem{theo}{Théorème}
\numberwithin{theo}{section}
\newtheorem{prop}[theo]{Proposition} 
\newtheorem{coro}[theo]{Corollaire}
 
\newtheorem{déf}[theo]{Définition}
\newtheorem{conje}[theo]{Conjecture} 
\newtheorem{question}[theo]{Question} 
\newtheorem{exmp}[theo]{Exemple} 
\newtheorem*{nota}{Notation}

 \newtheorem{rem}[theo]{Remarque} 

\newcommand{\R} {\ensuremath{\mathbb{R}}}
\newcommand{\Z} {\ensuremath{\mathbb{Z}}}
\newcommand{\h} {\ensuremath{\mathbb{H}}}
\newcommand{\F} {\ensuremath{\mathcal{F}}}

\newcommand{\borg} {\ensuremath{\partial \Gamma}}
\newcommand{\N} {\ensuremath{\mathbb{N}}}

\newcommand{\s} {\ensuremath{\mathbb{S}}}

\newcommand{\modcont}[2] {\ensuremath{\mathrm{Mod}_#1} (#2)}
\newcommand{\modcomb}[2] {\ensuremath{\mathrm{Mod}_#1} (#2)}

\newcommand{\modcombg}[1] {\ensuremath{\mathrm{Mod}_p} (#1,G_k)}

\newcommand{\modcombfo} {\ensuremath{\mathrm{Mod}_p} (\mathcal{F}_0,G_k)}

\newcommand{\di}[1] {\mathrm{dist}(\ensuremath{#1})}
\newcommand{\confdim}[1] {\mathrm{Confdim}(\ensuremath{#1})}

\newcommand{\dia}[1] {\mathrm{diam} \ensuremath{\:#1}}

   \newcommand{\numerote} [1] {
  \begin {enumerate}
  #1
  \end{enumerate}}
  
   \newcommand{\numeroti} [1] {
  \begin {enumerate} [i)]
  #1
  \end{enumerate}}
  
     \newcommand{\liste} [1] {
  \begin {itemize}#1
  \end{itemize}}

\title{Propriétés  combinatoires du  bord d'un groupe hyperbolique}
\author{Antoine Clais}
\address{Technion\\Department of Mathematics\\32000 Haifa, Israel} 
\email{aclais@tx.technion.ac.il}
\date{\today} 

\numberwithin{equation}{section}
  
\begin{document}
    \selectlanguage{french}
   \maketitle
 \begin{abstract}
Le but de ce survol est de présenter les modules combinatoires récemment utilisés pour étudier les propriétés quasi-conformes des bords des groupes hyperboliques.  Dans un premier temps, on rappellera quelques résultats et questions de rigidité bien connus qui ont motivés l'introduction de ces outils.  Puis on définira les modules combinatoires et la propriété de Loewner combinatoire qui offrent une nouvelle approche pour résoudre des problèmes ouverts depuis longtemps. Enfin, on décrira des applications concrètes de ces outils à travers quelques résultats récents et questions ouvertes.  
 \end{abstract}

   \selectlanguage{english}
 \begin{abstract}
 The goal of this survey is to present combinatorial modulus that have been recently used to study quasi-conformal properties of boundaries of hyperbolic groups. First, we will recall well known rigidity results and questions that motivated the introduction of these tools. Then we will define combinatorial modulus and the combinatorial Loewner property that provide new approaches to rigidity questions. Finally, we will describe some applications of these tools through recent results and open questions.
 \end{abstract}

   \selectlanguage{french}
\begin{flushleft}{\bf Mots-clés:} Bord d'un groupe hyperbolique, analyse quasi-conforme, modules combinatoires.\end{flushleft}
\begin{flushleft}{\bf 2010 Mathematics Subject Classification:}  	20F67, 30L10.\end{flushleft}

\setcounter{tocdepth}{1}
\tableofcontents

\section{Rigidité et homéomorphismes quasi-Möbius au bord} \label{section de QIaQM}
La géométrie d'un espace à courbure négative est intimement liée à la structure quasi-conforme de son bord. Ce fait a été utilisé dès les prémisses de la géométrie hyperbolique et de la théorie géométrique des groupes. En effet, en 1883 H. Poincaré, dans son \emph{Mémoire sur les groupes Kleinéens}, utilise le groupe conforme de la sphère $\s^2$ pour construire des exemples de groupes  agissant géométriquement sur l'espace hyperbolique réel $\h^3$. Le lien entre analyse quasi-conforme au bord et géométrie de l'espace hyperbolique a ensuite mené G.D. Mostow à son fameux théorème de rigidité. À la suite de G.D. Mostow, les travaux de  P. Pansu, J.W. Cannon,  M. Bourdon-H. Pajot  et  M. Bonk-B. Kleiner ont fait fructifier ce lien dans le cadre des groupes hyperboliques au sens de Gromov.  
 
 L'objet de cette partie est de rappeler le lien entre quasi-isométries et homéomorphismes quasi-Möbius du bord à travers quelques résultats et questions classiques. On peut se référer à \cite{KleinerAsymptoticGeom},  \cite{HaissinskyGeomQConf} ou \cite{BourdonMostowType} pour des survols complets  sur ces sujets.
 
\subsection{Questions de rigidité} La théorie des espaces et des groupes hyperboliques au sens de Gromov, fait l'objet d'excellentes introductions dans \cite{GhysHarpe} et \cite{CoorDelPapa}.
 
 Rappelons qu'un groupe de type fini   $\Gamma$ \emph{agit géométriquement} sur un espace métrique $X$, si :  $\Gamma < \text{Isom}(X)$, le quotient $X/\Gamma$ est compact et  l'action de $\Gamma$ est proprement discontinue. 
Dans l'article fondateur \cite{GromovHyperGroups}, M. Gromov introduit la notion d'espace hyperbolique qui généralise les variétés à courbures strictement négatives. Un espace métrique propre et géodésique $X$ est dit \emph{hyperbolique} (au sens de Gromov) si tout ses triangles géodésiques sont fins. C'est-à-dire s'il existe une constante $\delta \geq 0$ telle que pour tout triangle géodésique $ [x,y ] \cup  [y,z ] \cup   [z,x ] \subset X$, tout point $p\in  [x,y ]$ vérifie \[\di{p,  [y,z ] \cup   [z,x ]}\leq \delta.\]
Un groupe de type fini agissant géométriquement sur un espace hyperbolique est appelé un \emph{groupe hyperbolique}. Dans ce cas, ses graphes de Cayley sont des espaces métriques hyperboliques.  

Tout au long de cette partie,  on désigne par $(X,\vert \cdot -\cdot \vert)$ et $(Y,\vert \cdot -\cdot \vert)$ deux espaces métriques, géodésiques, propres et hyperboliques au sens de Gromov. 
Rappelons que deux applications $F,F' : X \longrightarrow Y$ sont à \emph{distance bornée} l'une de l'autre s'il existe une constante $K > 0$ telle que $\vert F(x) -F'(x)\vert$ pour tout $x\in X$. Dans ce cadre général, on peut se poser les questions suivantes.

\begin{question} \label{questionsde rigidité}

\text{ }

\numerote{\item  \textbf{Rigidité à la Mostow.} Un isomorphisme de groupe $\rho : \Gamma_1\longrightarrow \Gamma_2$ entre deux réseaux cocompactes de $\mathrm{Isom}(X)$ est-il toujours réalisé par  conjugaison par un élément de $\mathrm{Isom}(X)$ ? 

\item \textbf{Rigidité des quasi-isométries (Q.I).} Est-il vrai que toute quasi-isométrie $F:X\longrightarrow Y$ est à distance bornée d'une isométrie ?}
\end{question}

Il convient de remarquer que la rigidité quasi-isométrique est plus forte que la rigidité à la Mostow. Les exemples suivants  motivent les question précédentes. 

\begin{exmp}\label{exRigide}\numerote{\item Le célèbre théorème de G.D. Mostow affirme que pour $d\geq 3$ l'espace hyperbolique réel $\h^d$ vérifie le premier type de rigidité (cf. \cite{MostowRigidity}). 

 \item P. Pansu a montré que pour $d\geq 2$, l'espace hyperbolique quaternionien de dimension $d$ vérifie le second type de rigidité (cf. \cite{PansuMetriquesdeCCetQI}). 
\item Les immeubles fuchsiens vérifient le second type de rigidité (cf.  \cite{BourdonPajotRigi} pour le cas à angles droits  et  \cite{XieQIRigidity} pour le cas général).}
\end{exmp}

A la lumière du cas des espaces hyperboliques réels,  la rigidité n'est pas une propriété générique des espaces hyperboliques. La rigidité est renforcée quand la structure hyperbolique est enrichie par une autre structure : quaternioniènne  ou immobilière dans les exemples.

Le prochain exemple pourrait être obtenu en étudiant les immeubles hyperboliques à angles droits de dimension 3. Cette question est discutée dans la Partie \ref{secQuestions} où un bon candidat à la rigidité Q.I est donné.

\subsection{Des quasi-isométries aux homéomorphismes quasi-Möbius au bord} Expliquons maintenant la stratégie classique qui permet de résoudre ce genre de questions en passant par le bord. Rappelons que, par le lemme de   Švarc-Milnor,   les réseaux $\Gamma_1$ et $\Gamma_2$ de la Question \ref{questionsde rigidité} sont quasi-isométriques, en tant que groupes de présentations finis munis de la métrique du mot, à $(X,\vert \cdot -\cdot \vert)$. En conséquence,    les deux questions de rigidité se réduisent essentiellement à des questions sur les quasi-isométries. 

Des quasi-isométries on passe  au homéomorphismes du bord par le lemme de Morse. Fixons un point base $x_0\in X$. Le \emph{bord à l'infini} de $X$, noté $\partial X$ est l'ensemble des rayons géodésiques dans $X$ muni de la relation d'équivalence qui identifie deux rayons $r,r' : [0,+\infty[\longrightarrow X$ dès qu'il existe $K>0$ vérifiant $\vert r(t)-r'(t)\vert \leq K$ pour tout $t\in [0,+\infty[$.

 Le lemme de Morse affirme que l'image par une quasi-isométrie d'un rayon géodésique reste à distance bornée d'un rayon géodésique. Cela entraîne qu'une quasi-isométrie entre espaces hyperboliques s'étend continûment  en homéomorphisme au bord. En particulier, si $\Gamma$ est un groupe hyperbolique agissant sur $X$, son bord, c'est-à-dire le bord d'un de ses graphe de Cayley, s’identifie naturellement au bord de $X$ par la quasi-isométrie $: g\in \Gamma \longmapsto gx_0 \in X$.

 Les métriques visuelles des bords fournissent plus d’information sur  cet homéomorphisme. Pour voir cela on utilise les notions suivantes introduites dans \cite{VaisalaQMmaps}.
 
\begin{déf} Soient $Z$ et $Z'$ deux espaces métriques.
\numeroti{\item Pour $a,b,c,d \in Z$ distincts on définit   le \emph{birapport} suivant
\[ [a : b : c: d ]=\frac{d(a,b)}{d(a,c)}\cdot\frac{d(c,d)}{d(b,d)}.\] 
\item  Un homéomorphisme $f : Z \longrightarrow Z'$  est dit \emph{quasi-Möbius} ou \emph{Q.M} s'il existe un homéomorphisme $\phi :[0,+\infty[ \longrightarrow [0,+\infty[$ tel que pour tout  $a,b,c,d \in Z$ distincts \[[f(a):f(b):f(c):f(d)] \leq \phi( [a:b:c:d ]).\]Dans ce cas $f^{-1}$ est aussi un homéomorphisme Q.M et  $Z$ et $Z'$ sont dits \emph{quasi-Möbius équivalents} ou \emph{Q.M équivalents}.}  
\end{déf}
\begin{nota}
Tout au long de ce survol, pour deux  fonctions réelles $f$ et $g$,  on note  $f\asymp g$ s'il existe $A>0$ vérifiant \[A^{-1} f \leq g \leq A f.\]  
\end{nota}
  Rappelons que l'hyperbolicité implique que pour tout $a, b \in \partial X$ il existe une géodésique $(a, b)\subset X$ joignant $a$ à $b$. Cela permet, en particulier, de munir $\partial X$  d'une \emph{métrique visuelle}. C'est-à-dire une métrique $d$ telle qu'il existe $\alpha>1$  vérifiant pour tout $a,b\in \partial X$ 
\[   d(a,b)\asymp \alpha ^\ell ,\]
où $\ell=\di{x_0,(a, b)}$. Par ailleurs, la finesse des triangles dans $X$ entraîne qu'il existe une constante uniforme $K>0$ telle que  pour tout $n\geq 0$ assez grand
\[\ell- K \leq \vert a_n - x_0 \vert + \vert b_n - x_0 \vert - \vert a_n - b_n \vert \leq \ell +K,\]
où $a_n$ et $b_n$ sont des suites convergeant respectivement vers $a$ et $b$. Quitte à modifier $\alpha$,   on a donc
\begin{equation} \label{equametricvisu}
 d(a,b)\asymp \lim_{n\rightarrow\infty} \alpha^{(\vert a_n - x_0 \vert + \vert b_n - x_0 \vert - \vert a_n - b_n \vert )}.
\end{equation}

Choisissons maintenant $a_n, b_n, c_n $ et $d_n$ quatre suites dans $X$ convergeant respectivement vers $a,b,c,d \in \partial X$ distincts.  En utilisant la relation \ref{equametricvisu} on peut écrire \[ [a : b : c: d ] \asymp \lim_{n\rightarrow\infty} \alpha^{(\vert a_n - b_n \vert + \vert c_n - d_n \vert - \vert a_n - c_n \vert - \vert b_n - d_n \vert)}.\] En notant $L= \di{(a,c),(b,d)}$ 
 on obtient, quitte à modifier $\alpha$, \begin{equation}\label{equaBirap}
 [a : b : c: d ] \asymp \alpha^L.
\end{equation}

Maintenant considérons  $F :X \longrightarrow Y$ une $(\lambda,K)$-quasi-isométrie. On note   $f$ l’homéomorphisme de $\partial X$ vers $\partial Y$ induit par le lemme de Morse. Les inégalités  \[ \lambda^{-1}  \di{(a,c),(b,d)} -K \leq \di{(f(a),f(c)),(f(b),f(d))}\leq \lambda  \di{(a,c),(b,d)} +K \]ainsi que la relation \ref{equaBirap} entraînent  la proposition suivante.

\begin{prop}
Toute quasi-isométrie $F: X\longrightarrow Y$ se prolonge continûment en un homéomorphisme quasi-Möbius $f:\partial X\longrightarrow\partial Y$.
\end{prop}

L'action cocompacte du groupe $\Gamma$ permet de montrer que, réciproquement, tout homéomorphisme Q.M au bord est le prolongement d'une quasi-isométrie (cf. \cite{PaulinGroupeHyperParSonBord}).

 \subsection{Questions de rigidité au bord} On peut donc réduire l'essentiel des  problèmes de rigidité à la question suivante au bord : les homéomorphismes quasi-Möbius du bord sont-ils les prolongements d'isométries ?

La stratégie employée par P. Pansu puis M. Bourdon-H. Pajot pour leurs théorèmes de rigidité   est analogue à la stratégie de Mostow pour montrer la rigidité des $\h^d$ pour $d\geq 3$. Dans le cas de $\h^d$, où le bord est $\partial \h^d \simeq \s^{d-1}$, cette stratégie est la suivante.

\numerote{\item On montre un théorème de rigidité des homéomorphismes quasi-Möbius au bord. En utilisant les propriétés  dynamiques de l'action du  groupe au bord, on identifie les homomorphismes quasi-Möbius avec les applications conformes du bord.
\item On construit une isométrie à l'intérieur en utilisant le théorème de Liouville qui identifie le groupe conforme de la sphère $\s^{d-1}$  au groupe d'isométrie de $\h^d$ (cf. par exemple \cite[Section A.3.]{BenePetroHyperGeom}).}

Pour appliquer cette stratégie au bord  d'un espace hyperbolique au sens de Gromov, ce schéma  nécessite des adaptations techniques importantes. En particulier, le bord d'un tel espace n'est pas forcément une variété et la notion d'application conforme  au bord ne va pas de soi.   Cependant, comme cela est expliqué dans \cite[Partie 5.]{HaissinskyGeomQConf} la notion d'application conforme du bord se généralise  et  les  questions de rigidité se réduisent essentiellement à la question :  les homéomorphismes quasi-Möbius du bord sont-ils conformes ?

Les   \emph{espaces de Loewner} sont les espaces dans lesquels cette question trouve une réponse positive (cf. \cite[Théorèmes 2.15 et 5.11]{HaissinskyGeomQConf}).

 \begin{rem}  \label{remPvdeGromovThurston} Dans le cas de l'espace hyperbolique $\h^d$, $d\geq3$ une preuve de la rigidité de Mostow reposant sur le volume simplicial  a été donnée par M. Gromov puis généralisée par W. Thurston (cf. \cite[Chapter 6]{ThurstonNotesonthreeManifolds}). Cette preuve a  ensuite été « dualisée » dans \cite{BBIDualGromovThurstonMostow}  en utilisant la cohomologie bornée. Il serait intéressant de savoir si ces approches peuvent-être utilisées pour retrouver les résultats de rigidité des Exemples \ref{exRigide}.$ii)$ et $iii)$. Cela  fournirait une nouvelle approche pour étendre ces résultats.
  \end{rem}

\section{Espaces de Loewner}\label{sec EspdeLoewner}

Les modules des courbes sont des invariants naturels dans la théorie des applications quasi-conformes de $\R^d$ (cf. \cite{VaisalaLecture}). Les espaces de Loewner introduits par  Heinonen et Koskela dans \cite{HeinKoskQConf} constituent une généralisation de cette théorie dans le cadre abstrait des espaces métriques mesurés. On peut se référer à l'ouvrage \cite{HeinonenLect} pour une introduction complète sur ce sujet. 

En première approximation, on peut dire que les espaces de Loewner sont les espaces \emph{« contenant beaucoup de courbes rectifiables »}. Néanmoins, l'intuition  traduit difficilement cette notion très technique.

\subsection{Modules analytiques des courbes} Dans cette partie, $(Z,d,\mu)$ désigne un espace métrique mesuré . Dans la suite, on se focalisera sur le cas où $Z$ est le bord d'un espace hyperbolique. On suppose, de plus, que $Z$ est compact et \emph{$Q$-Ahlfors-régulier} ($Q$-AR ou AR)  pour $Q>1$. C'est-à-dire que pour tout $0<R\leq  \dia{Z}$ et toute boule  $B\subset Z$ de rayon $R$   \[  \mu  (B) \asymp R^Q .\] Remarquons que sous cette hypothèse, la mesure $\mu$ est comparable à la mesure de Hausdorff relative à la distance $d$. En particulier, $Q$ est la dimension de Hausdorff de $(Z,d)$.  
  
  On appelle \emph{courbe} de $Z$ une application continue $\gamma : [0,1] \longrightarrow Z$. On confond par la suite une courbe et son image.   Soit $\F$ un ensemble de courbes dans $Z$.  Une fonction borélienne $f : Z \longrightarrow [0,+\infty [$ est dite \emph{$\F$-admissible} si pour toute courbe rectifiable $\gamma \in \F$ \[ \int_{0}^1 f(\gamma(t))dt\geq 1.\] Remarquons que la notion d’admissibilité ne fait pas appelle à la mesure $\mu$ mais simplement à la   métrique de $Z$. 
 
\begin{déf} On appelle \emph{$Q$-module de $\F$}   \[ \modcont{Q}{\F}= \inf \Big\{ \int_X f^{Q} d\mu\Big\}\] 
où l'infimum est pris sur l'ensemble des fonctions $\F$-admissibles et avec la convention $\modcont{Q}{ \F }=0$ si $\F$ ne contient pas de courbe rectifiable.  
\end{déf}

La proposition suivante  nous permet de voir le module comme une \emph{« mesure extérieure »} sur l’ensemble des courbes de $Z$.

\begin{prop}[\cite{HeinKoskQConf}]   \text{   } 
\numerote{\item Pour deux ensembles de courbes  $\F_1\subset \F_2$ on a  \[\modcont{Q}{\F_1}\leq \modcont{Q}{\F_2}.\]
\item Pour des ensembles de courbes $\F_1, \dots , \F_n$ on a \[\modcont{Q}{\bigcup\limits_{i=1}^n \F_i}\leq \sum\limits_{i=1}^n \modcont{Q}{\F_i}.\]}
\end{prop}

 \subsection{Espace de Loewner} On appelle \emph{continuum} un sous-ensemble compact et connexe de $Z$. Un continuum est, de plus, dit \emph{non-dégénéré} s'il contient au moins deux points. Pour $A$ et $B$ deux continua non-dégénérés disjoints on note  $\F (A,B)$ l'ensemble de toutes les courbes qui joignent $A$ à $B$, c'est-à-dire   qui intersectent $A$ et $B$. On écrit aussi $\modcont{Q}{A,B}:=\modcont{Q}{\F (A,B)}$. Enfin on appelle \emph{distance relative entre $A$ et $B$}  
\[ \Delta(A,B)  =  \frac{\di{A,B}}{\min\{\dia{A},\dia{B} \} } .\]
Les espaces de Loewner sont définis de la manière suivante.
\begin{déf} \label{defLoewner}
On dit que  $(Z,d,\mu)$ est un \emph{$Q$-espace de Loewner} (ou vérifie la \emph{$Q$-propriété de Loewner}) s'il existe une fonction croissante   $\phi: ]0, + \infty[ \longrightarrow ]0, + \infty[ $ vérifiant pour toute paire de continua non-dégénérés disjoints $A$  et  $B$ de $Z$    \[ \phi (\Delta(A,B)^{-1}) \leq \modcont{Q}{ \F(A,B)}. \]
\end{déf}

S'il n'y a pas d’ambiguïté possible, on omet de préciser la dimension $Q$. Par ailleurs, on utilise parfois le terme  propriété de Loewner \emph{analytique} pour éviter la confusion avec la propriété de Loewner combinatoire définie dans la Partie \ref{subsecCLP}. On peut interpréter l'inégalité de la définition de la manière suivante : \begin{center}
\emph{« il y a toujours beaucoup de courbes rectifiables qui joignent deux continua »}.
\end{center}
Si $Z$  est un espace de Loewner, le comportement asymptotique de la fonction $\phi$ est donné par \cite[Theorem 3.6.]{HeinKoskQConf}. Pour $t$ suffisamment petit  $\phi(t)\approx \log\frac{1}{t}$, pour $t$ suffisamment grand $\phi(t)\approx (\log t)^{1-Q}$.
 
Par ailleurs, la structure d'espace $Q$-AR fournit une majoration des modules.  
\begin{theo}[{\cite[Lemma 3.14.]{HeinKoskQConf}}]
Il existe une constante $C>0$ vérifiant la propriété suivante. Soient $A$ et $B$ deux continua non-dégénérés disjoints. Soient $0<2r<R$ et $x\in Z$  tels que $A\subset\overline{B(x,r)}$ et $B\subset X\backslash B(x,R)$. Alors  
\[\modcont{Q}{A,B}\leq C\Big(\log\frac{R}{r}\Big)^{1-Q} .\]
\end{theo}

Cette seconde inégalité peut s’interpréter de la manière suivante : 
\begin{center}
\emph{« plus deux continua sont éloignés moins il y a de courbes rectifiables les joignant »}.
\end{center}

\subsection{Exemples et non-exemples} Les espaces suivants sont de Loewner.
\begin{exmp}
 \numeroti{
\item L'espace euclidien  $\R ^d$ pour $d\geq 2$. Ce résultat est due à C. Loewner pour $d\geq3 $ (cf. \cite{LoewnerConfCapa} ou \cite[Section 10.12.]{VaisalaLecture} pour une preuve simple). Ce premier exemple,  avec la caractérisation de la propriété de Loewner par les inégalités de Poincaré, implique que toute variété réelle compacte de dimension $d\geq 2$ est de Loewner (cf. \cite[Section 6.1.]{HeinKoskQConf}).
\item Les variété riemanniennes complètes non-compactes à courbure de Ricci positive et à croissance du volume optimale de dimension $d\geq 2$ (cf. \cite[Section 6.3.]{HeinKoskQConf}). 
\item Les bords  des immeubles fuchsiens (cf. \cite{BourdonPajotPoinc}).
\item Certains tapis non-auto-similaires, homéomorphes au tapis de Sierpiński (cf. \cite{MackayTysonTapis}).
}
\end{exmp}

Un célèbre théorème due à  B. Bowditch dit  que  le bords d'un groupe hyperbolique non-fuchsien possède des points de coupures locaux si et seulement si le groupe se  scinde sur un sous-groupe cyclique   ou en une extension HNN  (cf. \cite[Theorem 6.2.]{BowditchCutPoints}). La proposition suivante implique donc que les bords des groupes qui se scindent, ne sont pas de Loewner.
  
\begin{prop}[{\cite{HeinKoskQConf}}]\label{proppointdecoupure}
Si $Z$ vérifie  la propriété de Loewner, alors $Z$ n'a pas de point de coupure locale.  \end{prop} 
Cependant, dans \cite{BourdonPajotCohomoEspaceBesov} sont présentés  des exemples de groupes hyperboliques ne se scindant pas dont le bord n'est pas Q.M équivalent à un espace de Loewner.

\section{Dimension conforme} \label{secdimconf}
 Comme nous allons le voir, la Définition \ref{defLoewner} est difficile à vérifier au bord d'un espace hyperbolique car elle requiert la connaissance de la \emph{dimension conforme}.  La dimension conforme est un invariant quasi-Möbius introduit par P. Pansu dans \cite{PansuDimconf}. Ici $(Z,d)$ est un espace métrique compact muni de la mesure de Hausdorff associé à $d$. On désigne   $\dim_\mathcal{H} (Z,d)$    la dimension de Hausdorff de  $(Z,d)$.  

\subsection{Définitions} La \emph{jauge conforme Ahlfors-régulière} de $(Z,d)$ est l'ensemble des espaces métriques
\[\mathcal{J}_c(Z,d):= \{ (Z',\delta) : (Z',\delta) \text{ est AR et est Q.M équivalent à  }  (Z,d)\}. \] 
\begin{déf}
La \emph{dimension conforme  Ahlfors-régulière} (ou  \emph{dimension conforme}) de $(Z,d)$ est   
\[\confdim{Z,d} := \inf\{\dim_\mathcal{H}(Z',\delta) : (Z',\delta) \in   \mathcal{J}_c(Z,d)  \}  . \]
\end{déf}

De la même manière que la dimension topologique et la dimension de Hausdorff sont respectivement invariantes par homéomorphismes et par homéomorphismes bi-Lipschitz, la dimension conforme est invariante par homéomorphismes quasi-Möbius. Les inclusions entre ces trois ensembles d'applications induisent les inégalités suivantes
\[\dim_T (Z)\leq \confdim{Z,d}\leq \dim_\mathcal{H}(Z,d),\]
où $\dim_T(Z)$ désigne la dimension topologique de $Z$.

\begin{exmp} \label{exDimConf}
 On présente ici quelques résultats concernant la dimension conforme de bords de groupes hyperboliques.  
\numeroti{
\item Tout d'abord M. Bourdon a calculé la dimension conforme du bord d'un immeuble  fuchsien  à angles droits (cf. \cite{BourdonImHyperDimConfRigi}). Si   $\Delta(p,q)$  est l'immeuble dont le groupe de Coxeter est le groupe de réflexions d'un $p$-gone hyperbolique à angles droits et d'épaisseur $q\geq 2$ alors
 \[\confdim{\Delta(p,q)}= 1+ \frac{\log(q-1)}{\mathrm{Arg} \cosh \frac{p-2}{2}}.\] 
 
\item Un résultat de J. Mackay, à mettre en parallèle avec la Proposition \ref{proppointdecoupure}, dit que si le bord d'un groupe hyperbolique est connexe et n'a pas de point de coupure locale alors il   est de dimension conforme strictement plus grande que $1$  (cf.  \cite{MackayConfdimGreaterThanOne}).   

\item Enfin J. Mackay a calculé l’asymptotique de la dimension conforme pour les groupes aléatoires, dans le modèle à densité et  dans le modèle à peu de relateurs  (cf. \cite{MackayConfDimSmallCancRandGroups}). Dans ce dernier cas la dimension conforme  est $2$.}

\end{exmp}
 
\subsection{Dimension conforme et espaces de Loewner} La jauge et la dimension conforme jouent un rôle crucial dans un espace de Loewner comme le montre le résultat suivant due à J. Tyson.
 \begin{theo}[{\cite[Corollary 4.2.2.]{MackayTysonConfDim}}] \label{theo Tyson} Soient  $Q>1$ et  $Z$ un espace métrique mesuré $Q$-AR et $Q$-Loewner. Alors $\confdim{Z}=Q$. 
 \end{theo}

 Réciproquement, le bord d'un groupe  hyperbolique est de Loewner  si sa dimension conforme est atteinte dans sa jauge conforme. 
 \begin{theo}[\cite{BonkKleinerConfDimGromHypergrps}] Si $Z$ est un espace métrique mesuré $Q$-AR de dimension conforme $Q>1$ et si $Z$ est Q.M équivalent au bord d'un groupe hyperbolique, alors $Z$ est  $Q$-Loewner. 
 \end{theo}
 
 Pour étudier la propriété de Loewner au bord des groupe hyperboliques, il faut faire face à la difficulté que cette propriété n'est pas stable par homéomorphisme Q.M.

\begin{theo}[\cite{TysonQCandQS}]
Soient  $Z$ et $Z'$ deux espaces métriques compacts respectivement  $Q$-Loewner et $Q'$-AR. Supposons que $Z$ et $Z'$ soient Q.M équivalents. Alors $Z'$ est un espace de Loewner si et seulement si $Q=Q'$.
\end{theo}

Les homéomorphismes Q.M, en général, ne préservent pas la dimension de Hausdorff et donc la propriété de Loewner n'est pas stable par homéomorphisme Q.M.

La difficulté de calculer la dimension conforme ainsi que la non-stabilité de la propriété de Loewner par homéomorphismes Q.M rendent cette propriété difficile à établir au bord d'un groupe hyperbolique.

\subsection{La conjecture de Cannon} La conjecture suivante due à J.W. Cannon est une conjecture de géométrie des groupes qui contient le cas hyperbolique de la conjecture de géométrisation de Thurston.

\begin{conje}[{\cite[Conjecture 5.1.]{CanSwenCurv}}]\label{conj Cannon}
Soit $\Gamma$ un groupe hyperbolique dont le bord est homéomorphe à $\s^2$. Alors $\Gamma$ agit géométriquement sur $\h^3$.
\end{conje}
 
 Les travaux de G. Perelman et leurs généralisations ont montré la conjecture dans le cas où $\Gamma$ est le groupe fondamental d'une 3-variété mais le cas général reste ouvert.   Les progrès de  M. Bonk et B. Kleiner autour de cette conjecture ont fourni de nombreux outils et résultats d'analyse quasi-conforme au bord des espaces hyperboliques. En particulier, ils ont montré que la conjecture suivante est équivalente à celle de Cannon.

 \begin{conje}[Conjecture 1.2 \cite{BonkKleinerConfDimGromHypergrps}]
 Soit $\Gamma$ un groupe hyperbolique dont le bord est homéomorphe à $\s^2$. Alors la dimension conforme de $\borg$ est atteinte par une métrique de sa jauge conforme.
 \end{conje} 
 
 \begin{rem} La conjecture de Cannon a récemment été abordée par les sous-groupes de surface. Dans \cite{MarkovicCritCannConj}, il est montré que si $\Gamma$ est un groupe hyperbolique dont le bord est homéomorphe à $\s^2$ et   si $\Gamma$ contient \emph{« suffisamment »} de sous-groupes de surface quasi-convexes. Alors $\Gamma$ agit géométriquement sur $\h^3$.
 
   Dans \cite{BeekerLazarSphereBoundaries}, l'hypothèse de V. Markovic est affaiblie. Il est montré que si $\Gamma$ est  un groupe hyperbolique à un bout, dont la $\Z/2\Z$-cohomologie s'annule à l'infini  et si $\Gamma$  contient suffisamment de sous-groupes de surface quasi-convexes. Alors $\Gamma$ agit géométriquement sur $\h^3$. 
   \end{rem}

\section{Modules combinatoires et propriété de Loewner combinatoire} \label{section modcomb CLP}
Pour résumer les deux parties précédentes, on peut dire que le problème posé pour montrer la propriété de Loewner est de trouver une métrique qui réalise la dimension conforme. Cela permet de définir ensuite la \emph{ « bonne mesure »} avec laquelle travailler pour définir les modules des courbes. 

Pour contourner cette difficulté  les    modules combinatoires ont été introduits  par M. Bonk et B. Kleiner dans \cite{BonkKleinerQuasiSymParamofSpheres} en suivant des idées de P. Pansu et J.W. Cannon. L'idée est de remplacer la \emph{ « bonne mesure »}  par des  mesures discrètes. En prenant une suite de mesures discrètes de plus en plus fines, on espère trouver asymptotiquement une propriété qui ressemble à la propriété de Loewner. On peut se référer à \cite{BourdonKleinerCLP} pour une introduction complète de ces outils combinatoires.
 
Dans cette partie $(Z,d)$ est un espace métrique compact. Pour $z\in Z$ et $0<r\leq \dia{Z}$, on désigne par $B(z,r)$ la boule  ouverte de centre $z$ et de rayon $r$. En pratique, on veut appliquer les techniques suivantes au bord d'un groupe hyperbolique.

\subsection{Modules combinatoires}\label{subsecmodcomb} Il est à noter que les définitions suivantes sont rigoureusement parallèles à celles  de la Partie \ref{sec EspdeLoewner}. De plus, les propriétés élémentaires des modules combinatoires se démontrent de manières analogues à celles des modules analytiques.

  Pour $k\geq 0$  et  $\kappa > 1$, on appelle  \emph{$\kappa$-approximation de $Z$ à l'échelle $k$} un recouvrement fini de $Z$, noté  $G_k$, par des ouverts tels que pour tout  $ v\in G_k$ il existe  $z_v\in v$ vérifiant la propriété suivante :

\liste{ \item  $B(z_v,\kappa^{-1} 2^{-k}) \subset
v \subset B(z_v,\kappa 2^{-k}) $,
\item  pour tout  $v\neq w\in G_k$ on a $B(z_v,\kappa^{-1} 2^{-k})\cap B(z_w,\kappa^{-1} 2^{-k})= \emptyset$.}
Une suite de recouvrement  $\{G_k \}_{k\geq 0}$ est appelée une \emph{$\kappa$-approximation de $Z$}.

On fixe $\{G_k \}_{k\geq 0}$  une telle approximation et on construit des modules combinatoires, de courbes grâce à cette approximation.  Soient  $\rho: G_k \longrightarrow [0,+\infty[$ une fonction  positive et $\gamma$ une courbe dans $Z$.  La \emph{$\rho$-longueur} de   $\gamma$  est donnée par \[  L_\rho (\gamma) =  \sum_{\substack{ \gamma\cap v \neq \emptyset \\  v\in G_k}}\rho (v).\] Pour  $p\geq 1$, la \emph{$p$-masse} de $\rho$ est donnée par \[M_p(\rho) = \sum\limits_{v\in G_k} \rho (v)^p.\]
On fixe maintenant  $p\geq 1$. Soit  $\F$ un ensemble non-vide de courbes de $Z$. Une fonction $\rho: G_k \longrightarrow [0,+\infty[$ est dite \emph{$\F$-admissible} si $L_\rho(\gamma)\geq1$ pour toute courbe $\gamma\in \F$.
 
\begin{déf} On appelle \emph{$p$-module  $G_k$-combinatoire de $\F$}  \[ \modcombg{\F}= \inf \{ M_p(\rho)\}\] 
où  infimum est pris sur l'ensemble des fonctions $\F$-admissibles et avec la convention \[\modcombg{ \emptyset }=0.\] 
\end{déf}

Comme annoncé, les modules combinatoires sont des \emph{« mesures extérieures  discrètes sur les ensembles de courbes »}.  
\begin{prop}[{\cite[Proposition 2.1.]{BourdonKleinerCLP}}] \label{propmodbase}\text{   } 
\numerote{\item  Pour deux ensembles de courbes $\F_1\subset \F_2$ on a \[\modcombg{\F_1}\leq \modcombg{\F_2}.\]
\item Pour des ensembles de courbes  $\F_1, \dots , \F_n$  on a \[\modcombg{\bigcup\limits_{i=1}^n \F_i}\leq \sum\limits_{i=1}^n \modcombg{\F_i}.\]}
\end{prop}
Il est important de noter que, contrairement aux modules analytiques, les modules combinatoires ne distinguent pas les courbes rectifiables et non-rectifiables.

  À partir de maintenant, on suppose que $Z$ est   un espace métrique compact et \emph{doublant}. Cela signifie qu'il existe une constante uniforme $N\geq 1$ telle que toute boule de rayon $0< r \leq \dia{Z}$ puisse être recouverte par $N$ boules de rayon $r/2$. Dans un espace métrique doublant le comportement asymptotique du $p$-module $G_k$-combinatoire ne dépend pas du choix de l'approximation.  Comme cela est montré dans \cite[Proposition 3.3.]{BourdonKleinerCLP}, le bord d'un groupe hyperbolique muni d'une métrique visuelle est doublant. En pratique, on manipule donc les modules combinatoires à changement d'approximation près.  
  
\subsection{Propriété de Loewner combinatoire (CLP)}\label{subsecCLP} Dans ce paragraphe $(Z,d)$ est un espace métrique compact connexe par arc et doublant. On fixe  $\kappa >1$, $\{G_k\}_{k\geq 0}$ une $\kappa$-approximation de $Z$ et $p\geq 1$. La définition qui suit est l'analogue combinatoire de la Définition \ref{defLoewner}.

De nouveau, pour deux continua non-dégénérés et disjoints on note $\F (A,B)$ l'ensemble de toutes les courbes de $Z$ joignant  $A$ et $B$, et on écrit $\modcombg{A,B}:=\modcombg{\F (A,B)}$.  

\begin{déf}
Pour  $p>1$, on dit que  $(Z,d)$ vérifie la   \emph{$p$-Propriété de Loewner Combinatoire} (CLP) s'il existe deux fonctions croissantes $\phi$ et $\psi$ sur $]0, + \infty[$ avec  $\lim_{t\rightarrow 0} \psi (t) =0$, telles que  
\numeroti{ \item pour toute paire  de continua non-dégénérés et disjoints  $A$ et  $B$ et pour tout  $k\geq 0$ avec $2^{-k}\leq \min\{\dia{A},\dia{B} \}$ on a  
\[ \phi (\Delta(A,B)^{-1}) \leq \modcombg{ A,B}, \]
\item pour toute paire de boules ouvertes $B_1$, $B_2$ de même centre et telles que $B_1\subset B_2$, et pour tout $k\geq 0$ avec $2^{-k}\leq   \dia{B_1}$ on a  
\[  \modcombg{\overline{B_1},Z\backslash B_2} \leq \psi (\Delta(\overline{B_1},Z\backslash B_2)^{-1}). \]
}

\end{déf}

On peut interpréter les inégalités de la définition de manière analogue à la Partie \ref{sec EspdeLoewner}. 
 \begin{center}
\begin{minipage}[c]{13cm}  
i) \emph{« il y a toujours beaucoup de courbes joignant deux continua »},

ii)  \emph{« plus deux continua sont éloignés moins il y a de courbes les joignant »}.
\end{minipage}
\end{center}

La proposition suivante est à mettre en parallèle avec le Théorème \ref{theo Tyson}.

\begin{prop}[{\cite[Corollary 3.7.]{BourdonKleinerCLP}}]
Soit $(\borg, d)$ le bord d'un groupe hyperbolique muni d'une métrique visuelle. Si $\borg$ vérifie la $p$-propriété de Loewner combinatoire alors $p=\confdim{\borg}$.
\end{prop}
 
Comme la propriété de Loewner analytique,  la CLP est donc  une propriété qui apparaît lorsqu'on se place à la bonne dimension. Néanmoins, dans le cas analytique la métrique qui réalise la  dimension conforme sert aussi, en amont, à définir une mesure qui  permet de calculer les modules analytiques. Dans le cas combinatoire seul la métrique intervient pour définir les modules des courbes. Cela autorise, en particulier, la CLP à être stable par homéomorphismes Q.M.  
 
\begin{theo}[{\cite[Theorem 2.6.]{BourdonKleinerCLP}}]
\label{theo CLP QM}
 Si  $Z'$ est Q.M équivalent à un espace métrique compact  $Z$ vérifiant la $p$-propriété de Loewner combinatoire, alors $Z'$ vérifie  aussi la $p$-propriété de Loewner combinatoire.
\end{theo}

Cette invariance en fait une propriété plus naturelle et plus facile à vérifier au bord des groupes hyperboliques. Elle permet notamment  de formuler une condition à satisfaire  en toute dimension qui implique la CLP (voir Proposition \ref{prop critéreCLP}).
 
\subsection{Propriétés de Loewner combinatoire et analytique} Dans le cas où $Z$ est un espace métrique mesuré compact et AR, les modules combinatoires et analytiques sont essentiellement les mêmes.

\begin{prop}[{\cite[Prop B.2.]{HaissinskyEmpilCercles}}]
 Soit $Z$ un espace métrique mesuré compact et $Q$-AR. Supposons que  $Z$ soit muni d'une approximation $\{G_k\}_{k\geq 0}$. Alors pour toute paire de continua non-dégénérés et disjoints $A,B$  et pour  $k\geq 0$ assez grand on a  
\[ \modcomb{Q}{A,B,G_k}\asymp \modcont{Q}{A,B}  \]
dès que  $\modcont{Q}{A,B}>0$ et $\lim \modcomb{Q}{A,B,G_k} = 0$ sinon. 
\end{prop}

En particulier, la propriété de Loewner implique la propriété de Loewner combinatoire. 
\begin{theo}[{\cite[Theorem 2.6.]{BourdonKleinerCLP}}]\label{theo loewnervers clp}
  Si $Z$ est un espace métrique compact   $Q$-AR   et Loewner, alors  $Z$ vérifie la $Q$-propriété de Loewner combinatoire.
\end{theo}

La conjecture suivante justifie l'attention portée à la CLP au bord des espaces hyperboliques.
\begin{conje} [{\cite[Conjecture 7.5.]{KleinerAsymptoticGeom}}]\label{conj CLPloewner}
Supposons que  $Z$ soit  quasi-Möbius équivalent au bord d'un groupe hyperbolique. Si  $Z$ vérifie la CLP alors $Z$ est quasi-Möbius équivalent à un espace de Loewner.
\end{conje}
  
Si la conjecture est vraie, la CLP au bord de groupes hyperboliques entraînera de nouveaux résultats de rigidité. Les espaces suivants vérifient la CLP sans que l'on sache s'ils sont quasi-Möbius équivalents à des espaces de Loewner.
\begin{exmp} \label{ex CLP}
\text{ }
\numeroti{ 
\item Le tapis de Sierpiński, l'éponge de Menger et leurs généralisations en dimensions supérieures plongés dans l'espace euclidien (cf. \cite[Theorem 4.1.]{BourdonKleinerCLP}),
\item  des bords de groupes de Coxeter de plusieurs types : les groupes simpliciaux, des groupes prismaux, des groupes hautement symétriques et des groupes à bords plans (cf. \cite[Section 8.]{BourdonKleinerCLP}),
\item des bords d'immeubles hyperboliques à angles droits de dimension 3 et 4 homéomorphe à des objets universels de Menger (cf. \cite{ClaCLPonBRAHB}),
}

\end{exmp}

\begin{figure}[h!]
 
\centering
\includegraphics[scale=0.24]{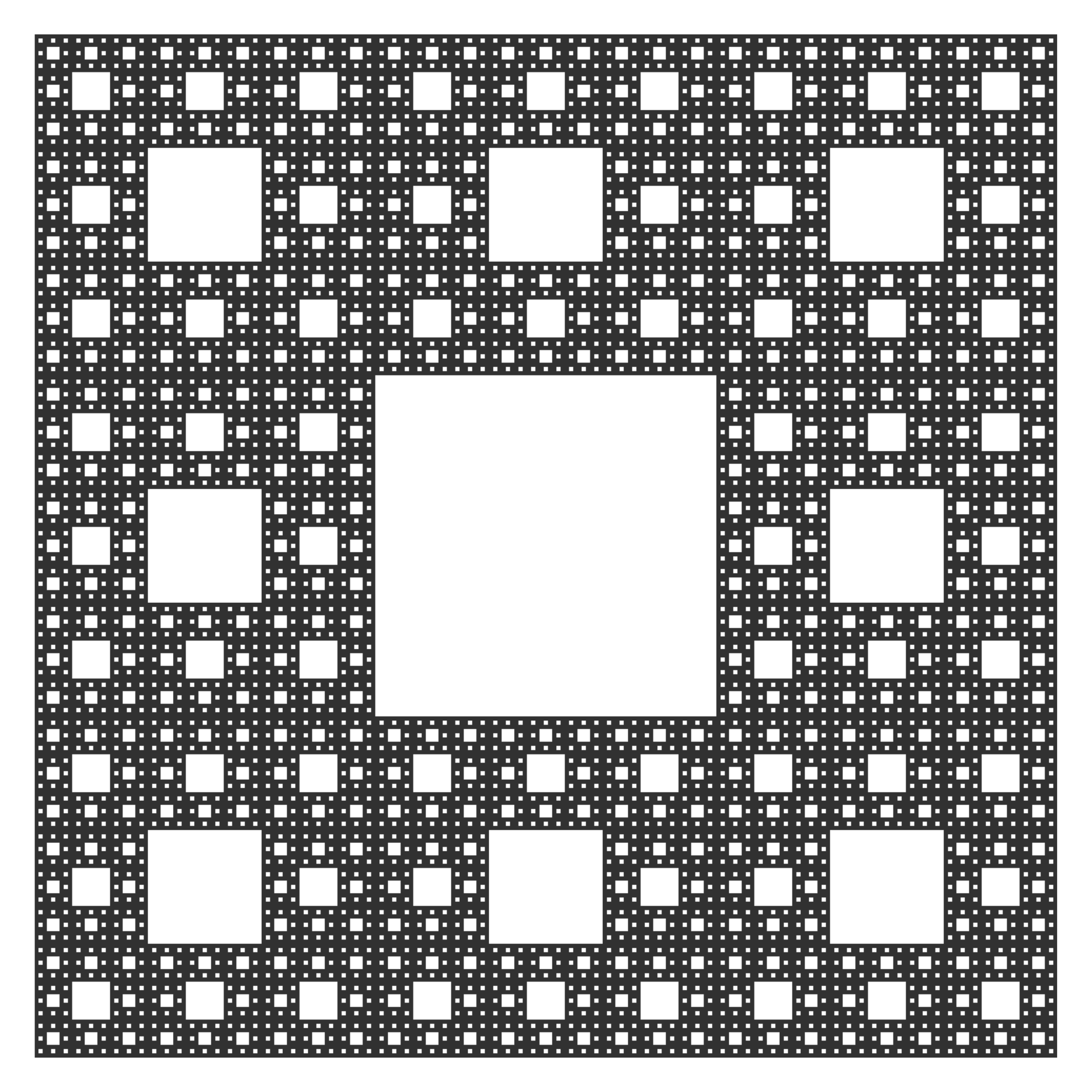}
\caption{Le tapis de  Sierpiński  vérifie la CLP. Est-il Q.M équivalent à un espace de Loewner ? Quelle-est sa dimension conforme ?}
\label{fig:D18}
\end{figure}

Toute résolution de la Conjecture \ref{conj CLPloewner} dans un cas particulier des exemples précédents serait un résultat intéressant.
Dans le cas des espaces auto-similaires,  le tapis de Sierpiński devrait être le plus facile à étudier. Pour l'instant, B. Kleiner  dans un travail non-publié et  indépendamment S. Keith et T. Laakso dans   \cite{KeithLaaksoConfAssDim} ont montré que la métrique euclidienne  sur le tapis de Sierpiński ne réalise pas la dimension conforme. La dimension  conforme du tapis reste inconnue.

 Les  immeubles hyperboliques  de l'Exemple \ref{ex CLP}.$iii)$ sont fortement suspectés d'être rigides (voir Partie \ref{secQuestions}). On pourrait prouver cela en montrant  que leurs bords sont de Loewner. 
  
\section{Utilisation des modules combinatoires aux bords} \label{secAppliModcomb}
Jusqu'à la fin de ce survol, $\Gamma$ un groupe hyperbolique non-élémentaire à bord connexe, $d$ désigne une métrique visuelle sur $\borg$,   $d_0>0$ est une constante suffisamment petite devant les constantes géométriques de $\borg$ et $\F_0$  est l'ensemble de toutes les courbes de $\borg$ de diamètre au moins $d_0$. Comme nous allons le voir, le comportement asymptotique des modules combinatoires se réduit au comportement asymptotique du module $\modcomb{p}{\F_0,G_k}$.

Plusieurs applications présentées ici concernent les groupes de Coxeter hyperboliques et les immeubles hyperboliques à angles droits.  On peut se référer à \cite[Section 5]{BourdonKleinerCLP} et à  \cite[Section 5]{ClaCLPonBRAHB} pour des introductions sur les groupes de Coxeter et les immeubles hyperboliques en vu de ces applications.

 Indépendamment de la Conjecture \ref{conj CLPloewner}, les modules combinatoires  sont  utiles en eux-mêmes pour étudier la structure quasi-conforme du bord d'un groupe hyperbolique. 
  
 \subsection{Bord sphérique} Dans \cite{BonkKleinerConfDimGromHypergrps} il est établi une condition sur les modules combinatoires pour que le bord d'un groupe homéomorphe à $\s^2$ soit Q.M équivalent à la sphère euclidienne. Combiné à \cite[Proposition 3.4.]{BourdonKleinerCLP} et au théorème p. 468 de \cite{SullivanErgoTheoInfiDiscGroupHyper}, cela  entraîne le corollaire suivant.
 
 \begin{coro}
  Soit $\Gamma$ un groupe hyperbolique dont le bord est homéomorphe à $\s^2$.  Supposons qu'il existe $C\geq 1$ tel que pour tout $k\geq 1$ on ait \[ \modcomb{2}{\F_0,G_k}\leq C. \]
  Alors $\Gamma$ agit géométriquement sur $\h^3$.
 \end{coro}
 
 Dans \cite{BourdonKleinerCLP}, ce corollaire permet de donner une nouvelle preuve de la conjecture de Cannon pour les groupes de Coxeter.
 
 \subsection{Dimension conforme} Par ailleurs,   S. Keith et B. Kleiner dans un travail non-publié, puis M. Carrasco dans \cite{CarrascoConfGauge}, dans un contexte plus général, ont montré  que la dimension conforme est égale à un exposant critique pour les modules combinatoires. Plus précisément,
 \begin{equation} \label{defEquivdimconf} \confdim{\borg} = \inf \{p\in [1, +\infty) : \lim_{k\rightarrow  + \infty} \modcombfo = 0\}.
  \end{equation}
 
 Dans le cas des immeubles hyperboliques à angles droits, la structure combinatoire au bord de l'immeuble est contrôlée par celle au bord d'un appartement. Cela implique que la dimension conforme du bord de l'immeuble est égale à un exposant critique calculé au bord d'un appartement (cf. \cite{ClaCLPonBRAHB}).
 
  Plus précisément, si $\Delta$ est un immeuble hyperbolique à angles droits de type $W$ et d'épaisseur constante $q\geq 2$, on note $G_k^W$  une approximation du bord $\partial W$ du bord du groupe de Coxeter,   $\mathrm{Mod}_p^W(\cdot)$ le  $p$-module  combinatoire calculé dans $\partial W$ et $\F_0^W$ l'ensemble de toutes les courbes de $\partial W$ de diamètre au moins $d_0$. Alors \[\confdim{\partial \Delta} = \inf \{p\in [1, +\infty) : \lim_{k\rightarrow  + \infty} (q-1)^k\mathrm{Mod}_p^W(\F_0^W,G_k^W) = 0\}.\]
  Dans \cite{ClaisConfDiamRAB} cette caractérisation est utilisée pour contrôler la dimension conforme du bord de l'immeuble par la dimension conforme du bord d'un appartement 
  
\begin{equation} \label{dimConImmeuble}
\confdim{\partial W} \cdot  \big( 1+\frac{\log(q-1)}{\tau}\big) \leq \confdim{\partial \Delta} \leq C \log(q-1),  
\end{equation} 
où    $\tau = \limsup_k \frac{1}{k} \log(\#\{g\in W : \vert g \vert \leq k \})$ est le taux de croissance de $W$ et $C>0$  une constante  indépendante de  $q$.  L'inégalité de gauche est optimale dans le cas fuchsien.

\section{Une méthode géométrique pour prouver la CLP} \label{sec methode geomCLP}
La proposition suivante permet d'établir la CLP sans avoir connaissance de la dimension conforme et de dégager un critère géométrique pour la CLP. 

\begin{prop}[{\cite[Proposition 4.5.]{BourdonKleinerCLP}}]\label{prop critéreCLP}
Soit $\{G_k\}_{k\geq0}$ une $\kappa$-approximation de $\borg$ muni d'une métrique visuelle.  Pour $p=1$, on suppose que $\modcombfo$ est non-borné. Pour $p\geq 1$, on suppose que pour toute courbe non-constante  $\eta \subset \borg$ et tout $\epsilon>0$, il existe   $C=C(p,\eta,\epsilon)$ tel que pour tout  $k \in \N$: \[ \modcombfo \leq C\cdot \modcomb{p}{\mathcal{U}_\epsilon(\eta),G_k}. \] 
Supposons de plus que $C$ puisse être choisi indépendamment de $p$ lorsque $p$ appartient à un compact de $[ 1, + \infty[$. Alors  $\borg$  vérifie la CLP.
\end{prop}

  Ici, pour une courbe $\eta$ de $\borg$, on désigne par  $\mathcal{U}_\epsilon (\eta)$ le $\epsilon$-voisinage de $\eta$ pour la distance $C^0$. C'est-à-dire qu'une courbe $\eta' \in \mathcal{U}_\epsilon(\eta)$ si et seulement si, il existe une paramétrisation $s:t \in [0,1] \longrightarrow [0,1]$ de $\eta$ telle que $d(\eta(s(t)),\eta'(t))<\epsilon$ pour tout  $t\in  [0,1]$.

  Dans le cas des immeubles hyperboliques à angles droits,   la structure combinatoire au bord de l'immeuble est contrôlée par celle au bord d'un appartement. En particulier,  \cite[Théorèmes 8.9 et 9.1]{ClaCLPonBRAHB} combinés à la proposition ci-dessus impliquent le corollaire suivant.
  
  \begin{coro}
  Soit $\Delta$ un immeuble hyperbolique à angles droits de type $W$ et d'épaisseur constante $q\geq 2$. Supposons que les hypothèses de la Proposition \ref{prop critéreCLP} soient vérifiées dans $\partial W$, alors $\partial \Delta$  vérifie la CLP.
  \end{coro}

Revenant au cas général, la première hypothèse de la Proposition \ref{prop critéreCLP} est vérifiée  si le bord est de dimension de Hausdorff strictement plus grande que $1$. En effet, dans ce cas on peut montrer  que pour tout  $N\geq1$  il  existe $N$ courbes disjointes de diamètre au moins $d_0$ dans $\borg$. Cela entraîne que pour $k\geq1$ suffisamment grand  $\modcomb{1}{\F_0,G_k}>N$.

Par les arguments habituellement utilisés pour comparer les modules combinatoires de deux ensembles de courbes (cf. par exemple la preuve de \cite[Theorem 6.1.]{BourdonKleinerCLP} ou de \cite[Theorem 6.12.]{ClaCLPonBRAHB}), la deuxième hypothèse   est vérifiée si la propriété suivante est vérifiée. 
 
\begin{center}\begin{minipage}[c]{12cm} 
$(S)  $ : Pour tout  $\epsilon>0$ et pour toute courbe $\eta$ de $\borg$, il existe $F$ un ensemble fini d'homéomorphismes bi-Lipschitz $f :\borg \longrightarrow \borg$, tel que pour toute courbe $\gamma\in \F_0$  le sous-ensemble $\bigcup_{f\in F} f(\gamma)$ de $\borg$ contient une courbe qui appartient à $\mathcal{U}_\epsilon (\eta)$.
\end{minipage}
\end{center}

Dans les groupes de Coxeter, les symétries de la chambre de Davis se prolongent au bord du groupe. Quand  ces symétries sont suffisamment nombreuses, on peut vérifier la propriété $(S)$ et montrer la CLP. 

\section{Questions} \label{secQuestions}
Cette dernière partie résume les problématiques abordées précédemment en insistant sur des exemples.

\subsection{Rigidité des immeubles de dimension supérieure}
Notons $D$ le dodécaèdre  régulier  à angles droits de $\h^3$, notons  $W_D$ le groupe de réflexions engendré par les faces de $D$ et notons enfin $\Delta_D$   un immeuble   de type $W_D$ et d'épaisseur constante $q\geq 2$.

\begin{question} \label{QuestionRigidité}
L'immeuble $\Delta_D$ est-il Mostow rigide ? Ou même Q.I rigide ?
\end{question}

Les réponses à ces questions ont de grandes chances d'être positives. En effet, les immeubles fuchsiens (les analogues de $\Delta_D$ en dimension 2) vérifient la rigidité Q.I alors même que leurs appartements ne sont pas rigides, ces derniers étant des plans hyperboliques $\h^2$. Les appartements de   $\Delta_D$ sont au contraire Mostow rigides, ce sont des espaces hyperboliques $\h^3$. La rigidité des appartements devrait donc être renforcée par la structure immobilière.

On sait déjà que le bord de $\Delta_D$ vérifie la CLP (voir Exemple \ref{ex CLP}.$iii)$). Cependant, la résolution de la Conjecture \ref{conj CLPloewner} dans ce cas particulier s'annonce difficile. Une étape décisive serait de trouver une métrique visuelle qui réalise la dimension conforme. Dans le cas fuchsien, il est fait usage de manière cruciale du fait suivant : pour obtenir une métrique sur un cercle, il suffit de disposer  d'une mesure sur ce cercle, et de définir la distance entre deux points comme la mesure  de l'arc de cercle qui les joint (cf. \cite[Lemme 3.1.4]{BourdonImHyperDimConfRigi}). Ce fait ne peut bien sûr pas être étendu à la sphère de dimension 2.

Une étape vers la résolution de ce problème serait d'améliorer les inégalités \ref{dimConImmeuble}. L'inégalité de gauche est optimale en dimension 2 mais il est difficile d’intuiter si elle sera optimale dans le cas de $\Delta_D$. Quoiqu'il en soit, l'amélioration des ces inégalités reposerait sur une compréhension très fine du comportement des modules combinatoires au bord d'un appartement. Par ailleurs, comme mentionné dans la Remarque \ref{remPvdeGromovThurston},  une approche reposant sur les volumes simpliciaux est aussi envisageable. 
 
Notons enfin que la question et les remarques précédentes trouvent un analogue  en dimension $4$ en substituant l'hécatonicosachore  (aussi appelé le 120-cellules) régulier  à angles droits de $\h^4$ au dodécaèdre. Pour l'instant, il est raisonnable de penser que la résolution de la Question \ref{QuestionRigidité} pour les immeubles « dodécaèdriques » permettra sa résolution pour les immeubles « hécatonicosachoriques ».

\subsection{Dimension conforme du bord d'un immeuble ni épais ni fin} Nous appelons maintenant immeuble d'épaisseur \emph{intermédiaire} un immeuble de type $(W,S)$ dont toutes les cloisons d'un type donné $s\in S$ sont d'épaisseurs $q\geq 3$ et toutes les cloison de type $s'\in S\backslash \{s\}$ sont fines. Le groupe de Coxeter du dodécaèdre  régulier  à angles droits de $\h^3$ est toujours noté $W_D$ et   $\Delta_{int}$ désigne  un immeuble de type $W_D$ d'épaisseur intermédiaire.

\begin{question} La dimension conforme de $\partial \Delta_{int}$ est-elle 2 ?\end{question}

De nouveau, cette question est suggérée par le cas fuchsien. Dans ce cas, le bord d'un immeuble d’épaisseur intermédiaire possède des points de coupures locaux et sa dimension conforme est 1 (voir Exemple \ref{exDimConf}.$ii)$).  En dimension 3, les points de coupures locaux sont remplacés par des « cercles de coupures locaux ». De manière générale, il est beaucoup plus difficile de comprendre l'influence d'un cercle de coupure locale sur le groupe. Cependant, dans ce cas précis  une analyse fine du module combinatoire au bord de l’appartement pourrait servir à répondre positivement à la question. Cela  constituerait un exemple intéressant de groupe avec un bord non-sphérique de dimension conforme 2.

 De nouveau, on peut substituer l'hécatonicosachore au dodécaèdre et se demander si le bord de l'immeuble d'épaisseur intermédiaire ainsi obtenu est de dimension conforme 3.

\subsection{Généricité} Disons pour conclure que la CLP et, a fortiori, la propriété analytique de Loewner sont des propriétés qui apparaissent de manière exotique au bord des groupes. Pour le moment, nous ne connaissons pas de cadre dans lequel ces propriétés soient génériques. Tout résultat de ce type serait un progrès majeur pour l'utilisation de ces techniques. Cependant rien n'indique qu'un  tel cadre existe et  tous les nouveaux exemples, aussi exotiques soient ils, sont intéressants à étudier.
 
 \begin{question}
Dans quelle classe de groupes la CLP est-elle une propriété générique ?
\end{question}
 
Des candidats naturels pourraient être les groupes aléatoires, dans le modèle à peu de relateurs ou à densité, dont la dimension conforme est déjà connue.

\begin{rem} La topologie au bord des immeubles d’épaisseur intermédiaire  $\Delta_{int}$   est plus rigide que celle au bord des immeubles épais. En particulier, elle suffit à montrer que ces immeubles vérifient la rigidité de Mostow (cf. \cite{LafontRigiditySurvey}).
\end{rem}

\subsection*{Remerciements} L'auteur remercie vivement Pierre Will,  pour lui avoir donné l’opportunité de publier ce survol, ainsi que Marc Bourdon, pour les quatre années de discussions qui ont nourri ce texte. Il est aussi reconnaissant au relecteur pour ses remarques et corrections précises.

   \bibliography{Biblio}

\begin{thebibliography}{MTW13}

\bibitem[BBI13]{BBIDualGromovThurstonMostow}
Michelle Bucher, Marc Burger, and Alessandra Iozzi.
\newblock A dual interpretation of the {G}romov-{T}hurston proof of {M}ostow
  rigidity and volume rigidity for representations of hyperbolic lattices.
\newblock In {\em Trends in harmonic analysis}, volume~3 of {\em Springer INdAM
  Ser.}, pages 47--76. Springer, Milan, 2013.

\bibitem[BK02]{BonkKleinerQuasiSymParamofSpheres}
Mario Bonk and Bruce Kleiner.
\newblock Quasisymmetric parametrizations of two-dimensional metric spheres.
\newblock {\em Invent. Math.}, 150(1):127--183, 2002.

\bibitem[BK05]{BonkKleinerConfDimGromHypergrps}
Mario Bonk and Bruce Kleiner.
\newblock Conformal dimension and {G}romov hyperbolic groups with 2-sphere
  boundary.
\newblock {\em Geom. Topol.}, 9:219--246, 2005.

\bibitem[BK13]{BourdonKleinerCLP}
Marc Bourdon and Bruce Kleiner.
\newblock Combinatorial modulus, the combinatorial {L}oewner property, and
  {C}oxeter groups.
\newblock {\em Groups Geom. Dyn.}, 7(1):39--107, 2013.

\bibitem[BL16]{BeekerLazarSphereBoundaries}
Benjamin Beeker and Nir Lazarovich.
\newblock Sphere boundaries of hyperbolic groups.
\newblock {\em Pre-print arXiv:1512.00866}, 2016.

\bibitem[Bou]{BourdonMostowType}
Marc Bourdon.
\newblock Mostow type rigidity theorems.
\newblock {\em to appear in Handbook of Group Actions}, Vol. III.

\bibitem[Bou97]{BourdonImHyperDimConfRigi}
Marc Bourdon.
\newblock Immeubles hyperboliques, dimension conforme et rigidit\'e de
  {M}ostow.
\newblock {\em Geom. Funct. Anal.}, 7(2):245--268, 1997.

\bibitem[Bow98]{BowditchCutPoints}
Brian~H. Bowditch.
\newblock Cut points and canonical splittings of hyperbolic groups.
\newblock {\em Acta Math.}, 180(2):145--186, 1998.

\bibitem[BP92]{BenePetroHyperGeom}
Riccardo Benedetti and Carlo Petronio.
\newblock {\em Lectures on hyperbolic geometry}.
\newblock Universitext. Springer-Verlag, Berlin, 1992.

\bibitem[BP99]{BourdonPajotPoinc}
Marc Bourdon and Herv{\'e} Pajot.
\newblock Poincar\'e inequalities and quasiconformal structure on the boundary
  of some hyperbolic buildings.
\newblock {\em Proc. Amer. Math. Soc.}, 127(8):2315--2324, 1999.

\bibitem[BP00]{BourdonPajotRigi}
Marc Bourdon and Herv{\'e} Pajot.
\newblock Rigidity of quasi-isometries for some hyperbolic buildings.
\newblock {\em Comment. Math. Helv.}, 75(4):701--736, 2000.

\bibitem[BP03]{BourdonPajotCohomoEspaceBesov}
Marc Bourdon and Herv{\'e} Pajot.
\newblock Cohomologie {$l_p$} et espaces de {B}esov.
\newblock {\em J. Reine Angew. Math.}, 558:85--108, 2003.

\bibitem[CDP90]{CoorDelPapa}
Michel Coornaert, Thomas Delzant, and Athanase Papadopoulos.
\newblock {\em G\'eom\'etrie et th\'eorie des groupes}, volume 1441 of {\em
  Lecture Notes in Mathematics}.
\newblock Springer-Verlag, Berlin, 1990.
\newblock Les groupes hyperboliques de Gromov. [Gromov hyperbolic groups], With
  an English summary.

\bibitem[Cla16a]{ClaCLPonBRAHB}
Antoine Clais.
\newblock Combinatorial {M}odulus on {B}oundary of {R}ight-{A}ngled
  {H}yperbolic {B}uildings.
\newblock {\em Anal. Geom. Metr. Spaces}, 4:Art. 1, 2016.

\bibitem[Cla16b]{ClaisConfDiamRAB}
Antoine Clais.
\newblock Conformal dimension on boundary of right-angled hyperbolic buildings.
\newblock {\em Pre-print arXiv:1602.08611}, 2016.

\bibitem[CP13]{CarrascoConfGauge}
Matias Carrasco~Piaggio.
\newblock On the conformal gauge of a compact metric space.
\newblock {\em Ann. Sci. \'Ec. Norm. Sup\'er. (4)}, 46(3):495--548 (2013),
  2013.

\bibitem[CS98]{CanSwenCurv}
James~W. Cannon and Eric~L. Swenson.
\newblock Recognizing constant curvature discrete groups in dimension {$3$}.
\newblock {\em Trans. Amer. Math. Soc.}, 350(2):809--849, 1998.

\bibitem[GdlH90]{GhysHarpe}
{\'E}tienne Ghys and Pierre de~la Harpe.
\newblock Espaces m\'etriques hyperboliques.
\newblock In {\em Sur les groupes hyperboliques d'apr\`es {M}ikhael {G}romov
  ({B}ern, 1988)}, volume~83 of {\em Progr. Math.}, pages 27--45. Birkh\"auser
  Boston, Boston, MA, 1990.

\bibitem[Gro87]{GromovHyperGroups}
M.~Gromov.
\newblock Hyperbolic groups.
\newblock In {\em Essays in group theory}, volume~8 of {\em Math. Sci. Res.
  Inst. Publ.}, pages 75--263. Springer, New York, 1987.

\bibitem[Ha{\"{\i}}09a]{HaissinskyEmpilCercles}
Peter Ha{\"{\i}}ssinsky.
\newblock Empilements de cercles et modules combinatoires.
\newblock {\em Ann. Inst. Fourier (Grenoble)}, 59(6):2175--2222, 2009.

\bibitem[Ha{\"{\i}}09b]{HaissinskyGeomQConf}
Peter Ha{\"{\i}}ssinsky.
\newblock G\'eom\'etrie quasiconforme, analyse au bord des espaces m\'etriques
  hyperboliques et rigidit\'es [d'apr\`es {M}ostow, {P}ansu, {B}ourdon,
  {P}ajot, {B}onk, {K}leiner{$\ldots$}].
\newblock {\em Ast\'erisque}, (326):Exp. No. 993, ix, 321--362 (2010), 2009.
\newblock S{\'e}minaire Bourbaki. Vol. 2007/2008.

\bibitem[Hei01]{HeinonenLect}
Juha Heinonen.
\newblock {\em Lectures on analysis on metric spaces}.
\newblock Universitext. Springer-Verlag, New York, 2001.

\bibitem[HK98]{HeinKoskQConf}
Juha Heinonen and Pekka Koskela.
\newblock Quasiconformal maps in metric spaces with controlled geometry.
\newblock {\em Acta Math.}, 181(1):1--61, 1998.

\bibitem[KL04]{KeithLaaksoConfAssDim}
Stephen Keith and Tomi~J. Laakso.
\newblock Conformal {A}ssouad dimension and modulus.
\newblock {\em Geom. Funct. Anal.}, 14(6):1278--1321, 2004.

\bibitem[Kle06]{KleinerAsymptoticGeom}
Bruce Kleiner.
\newblock The asymptotic geometry of negatively curved spaces: uniformization,
  geometrization and rigidity.
\newblock In {\em International {C}ongress of {M}athematicians. {V}ol. {II}},
  pages 743--768. Eur. Math. Soc., Z\"urich, 2006.

\bibitem[Laf07]{LafontRigiditySurvey}
Jean-Fran{\c{c}}ois Lafont.
\newblock Rigidity of hyperbolic {$P$}-manifolds: a survey.
\newblock {\em Geom. Dedicata}, 124:143--152, 2007.

\bibitem[Loe59]{LoewnerConfCapa}
Charles Loewner.
\newblock On the conformal capacity in space.
\newblock {\em J. Math. Mech.}, 8:411--414, 1959.

\bibitem[Mac10]{MackayConfdimGreaterThanOne}
John~M. Mackay.
\newblock Spaces and groups with conformal dimension greater than one.
\newblock {\em Duke Math. J.}, 153(2):211--227, 2010.

\bibitem[Mac14]{MackayConfDimSmallCancRandGroups}
John~M. Mackay.
\newblock Conformal dimension via subcomplexes for small cancellation and
  random groups.
\newblock {\em Math. Annalen. to appear. arxiv:1409.0802}, 2014.

\bibitem[Mar13]{MarkovicCritCannConj}
Vladimir Markovic.
\newblock Criterion for {C}annon's conjecture.
\newblock {\em Geom. Funct. Anal.}, 23(3):1035--1061, 2013.

\bibitem[Mos68]{MostowRigidity}
George~D. Mostow.
\newblock Quasi-conformal mappings in {$n$}-space and the rigidity of
  hyperbolic space forms.
\newblock {\em Inst. Hautes \'Etudes Sci. Publ. Math.}, (34):53--104, 1968.

\bibitem[MT10]{MackayTysonConfDim}
John~M. Mackay and Jeremy~T. Tyson.
\newblock {\em Conformal dimension}, volume~54 of {\em University Lecture
  Series}.
\newblock American Mathematical Society, Providence, RI, 2010.
\newblock Theory and application.

\bibitem[MTW13]{MackayTysonTapis}
John~M. Mackay, Jeremy~T. Tyson, and Kevin Wildrick.
\newblock Modulus and {P}oincar\'e inequalities on non-self-similar
  {S}ierpi\'nski carpets.
\newblock {\em Geom. Funct. Anal.}, 23(3):985--1034, 2013.

\bibitem[Pan89a]{PansuDimconf}
Pierre Pansu.
\newblock Dimension conforme et sph\`ere \`a l'infini des vari\'et\'es \`a
  courbure n\'egative.
\newblock {\em Ann. Acad. Sci. Fenn. Ser. A I Math.}, 14(2):177--212, 1989.

\bibitem[Pan89b]{PansuMetriquesdeCCetQI}
Pierre Pansu.
\newblock M\'etriques de {C}arnot-{C}arath\'eodory et quasiisom\'etries des
  espaces sym\'etriques de rang un.
\newblock {\em Ann. of Math. (2)}, 129(1):1--60, 1989.

\bibitem[Pau96]{PaulinGroupeHyperParSonBord}
Fr{\'e}d{\'e}ric Paulin.
\newblock Un groupe hyperbolique est d\'etermin\'e par son bord.
\newblock {\em J. London Math. Soc. (2)}, 54(1):50--74, 1996.

\bibitem[Sul81]{SullivanErgoTheoInfiDiscGroupHyper}
Dennis Sullivan.
\newblock On the ergodic theory at infinity of an arbitrary discrete group of
  hyperbolic motions.
\newblock In {\em Riemann surfaces and related topics: {P}roceedings of the
  1978 {S}tony {B}rook {C}onference ({S}tate {U}niv. {N}ew {Y}ork, {S}tony
  {B}rook, {N}.{Y}., 1978)}, volume~97 of {\em Ann. of Math. Stud.}, pages
  465--496. Princeton Univ. Press, Princeton, N.J., 1981.

\bibitem[Thu80]{ThurstonNotesonthreeManifolds}
William~P. Thurston.
\newblock {\em The Geometry and Topology of Three-Manifolds}.
\newblock Notes of Princeton University. 1980.

\bibitem[Tys98]{TysonQCandQS}
Jeremy~T. Tyson.
\newblock Quasiconformality and quasisymmetry in metric measure spaces.
\newblock {\em Ann. Acad. Sci. Fenn. Math.}, 23(2):525--548, 1998.

\bibitem[V{\"a}i71]{VaisalaLecture}
Jussi V{\"a}is{\"a}l{\"a}.
\newblock {\em Lectures on {$n$}-dimensional quasiconformal mappings}.
\newblock Lecture Notes in Mathematics, Vol. 229. Springer-Verlag, Berlin-New
  York, 1971.

\bibitem[V{\"a}i85]{VaisalaQMmaps}
Jussi V{\"a}is{\"a}l{\"a}.
\newblock Quasi-{M}\"obius maps.
\newblock {\em J. Analyse Math.}, 44:218--234, 1984/85.

\bibitem[Xie06]{XieQIRigidity}
Xiangdong Xie.
\newblock Quasi-isometric rigidity of {F}uchsian buildings.
\newblock {\em Topology}, 45(1):101--169, 2006.

\end{thebibliography}
  \bibliographystyle{alpha}

\end{document}